\magnification=1200


\font\tenbfit=cmbxti10 
\font\sevenbfit=cmbxti10 at 7pt 
\font\sixbfit=cmbxti5 at 6pt 

\newfam\mathboldit 

\textfont\mathboldit=\tenbfit
  \scriptfont\mathboldit=\sevenbfit
   \scriptscriptfont\mathboldit=\sixbfit

\def\bfit           
{\tenbfit           
   \fam\mathboldit 
}

\def\K{{\bf {K}}}

\def\Z{{\bf Z}}     
\def\R{{\bf R}}

\def\Bad{{\bfit Bad}}
\def\oU{{\overline U}}
\def\oV{{\overline V}}
\def\oW{{\overline W}}

\def\oA{{\overline A}}
\def\oB{{\overline B}}

\def\house#1{\setbox1=\hbox{$\,#1\,$}%
\dimen1=\ht1 \advance\dimen1 by 2pt \dimen2=\dp1 \advance\dimen2 by 2pt
\setbox1=\hbox{\vrule height\dimen1 depth\dimen2\box1\vrule}%
\setbox1=\vbox{\hrule\box1}%
\advance\dimen1 by .4pt \ht1=\dimen1
\advance\dimen2 by .4pt \dp1=\dimen2 \box1\relax}

  \def\eps{{\varepsilon}}

\def\sm{\smallskip} \def\ens{\enspace} 

\def\build#1_#2^#3{\mathrel{\mathop{\kern 0pt#1}\limits_{#2}^{#3}}}

\def\date {le\ {\the\day}\ \ifcase\month\or janvier
\or fevrier\or mars\or avril\or mai\or juin\or juillet\or
ao\^ut\or septembre\or octobre\or novembre
\or d\'ecembre\fi\ {\oldstyle\the\year}}

\font\fivegoth=eufm5 \font\sevengoth=eufm7 \font\tengoth=eufm10

\newfam\gothfam \scriptscriptfont\gothfam=\fivegoth
\textfont\gothfam=\tengoth \scriptfont\gothfam=\sevengoth

\def\smallsquare{\vbox{\hrule\hbox{\vrule height 1 ex\kern 1 ex\vrule}\hrule}}
\def\cqfd{\hfill \smallsquare\vskip 3mm}


\centerline{}

\vskip 4mm

\centerline{
\bf On the Littlewood conjecture in simultaneous
Diophantine approximation}

\vskip 8mm
\centerline{Boris A{\sevenrm DAMCZEWSKI} \footnote{}{\rm 
} \
\& \ Yann B{\sevenrm UGEAUD} 
\footnote{*}{Supported by the Austrian Science
Fundation FWF, grant M822-N12. } 
}

{\narrower\narrower
\vskip 12mm

\proclaim Abstract. {For any given real number $\alpha$ with
bounded partial quotients, we construct explicitly continuum many
real numbers $\beta$ with bounded partial quotients
for which the pair $(\alpha, \beta)$ satisfies
a strong form of the Littlewood conjecture. Our proof is 
elementary and rests on the basic theory of continued fractions.
}

}

\vskip 6mm
\vskip 8mm

\centerline{\bf 1. Introduction}

\vskip 6mm

It follows from the theory of continued fractions that, for 
any real number $\alpha$, there exist infinitely many positive integers $q$ 
such that
$$
q \cdot \Vert q \alpha \Vert < 1,  \eqno (1.1)
$$
where $\Vert \cdot \Vert$ denotes the distance to the nearest
integer. In particular, 
for any given pair $(\alpha,\beta)$ of real
numbers,  there exist infinitely many positive integers $q$ 
such that
$$
q \cdot \Vert q \alpha \Vert \cdot \Vert q\beta\Vert< 1.
$$
A famous open problem in simultaneous 
Diophantine approximation, called 
the Littlewood conjecture [8], claims that in fact, 
for any given pair $(\alpha, \beta)$ of real numbers, 
a stronger result holds, namely
$$
\inf_{q \ge 1} \, q \cdot \Vert q \alpha \Vert \cdot \Vert q \beta \Vert = 0.
\eqno (1.2)
$$

Throughout the present Note, we denote by $\Bad$ the set of badly
approximable numbers, that is,
$$
\Bad = \{ \alpha \in \R : \inf_{q \ge 1} \,
q \cdot \Vert q \alpha \Vert > 0\}.
$$
It is well-known that a real number lies in  
$\Bad$ if, and only if, it has bounded
partial quotients in its continued fraction expansion. 
It then follows that 
the Littlewood conjecture holds true for the pair $(\alpha, \beta)$
if $\alpha$ or $\beta$ has unbounded
partial quotients in its continued fraction expansion. It also
holds when the numbers $1$, $\alpha$, and $\beta$ are linearly dependent
over the rational integers, as follows from (1.1).

The first significant
contribution towards the Littlewood conjecture
goes back to Cassels \& Swinnerton-Dyer [3] who showed that
(1.2) holds when $\alpha$ and $\beta$ belong to the same cubic field.
However, since it is still not known
whether or not cubic real numbers have bounded
partial quotients, their result does not
yield examples of pairs of badly approximable real numbers
for which the Littlewood conjecture holds. 

In view of the above discussion, it is natural to restrict our
attention to independent parameters $\alpha$ and $\beta$,  
both lying in $\Bad$. 
The present paper is mainly devoted to the study of 
the following problem:

\medskip

{\it Question 1.} Given $\alpha$ in $\Bad$, is there any independent
$\beta$ in $\Bad$ so that the Littlewood conjecture is true for 
the pair $(\alpha, \beta)$?

\medskip

\noindent 
Apparently, Question 1 remained unsolved until 2000. It has then been
answered positively by Pollington \& Velani [12], who
established the following stronger result.

\medskip

\proclaim Theorem PV. Given $\alpha$ in $\Bad$, there exists a subset
$A(\alpha)$ of $\Bad$ with Hausdorff dimension one, such that,
for any $\beta$ in $A(\alpha)$, there exist infinitely many
positive integers $q$ with
$$
q \cdot \Vert q \alpha \Vert \cdot \Vert q \beta \Vert 
\le {1 \over \log q}.
\eqno (1.3)
$$
In particular, the Littlewood conjecture holds for 
the pair $(\alpha, \beta)$ for any $\beta$ in $A(\alpha)$.

\medskip

\noindent 
The proof of Theorem PV depends on sophisticated 
tools from metric number theory. At the end of [12],
Pollington \& Velani gave an alternative proof of a weaker
version of Theorem PV, namely with (1.3) replaced by (1.2). However,
even for establishing this weaker version, deep tools from
metric number theory are still needed, including a 
result of Davenport, Erd\H os and LeVeque on uniform distribution [4]  
and the {\it Kaufman measure} constructed in [7]. 

Very recently,
Einsiedler, Katok \& Lindenstrauss [6] proved the
outstanding result that the set of pairs of real numbers for which
the Littlewood conjecture does not hold has Hausdorff dimension zero.
Obviously, this implies a positive answer to Question 1. 
Actually, the authors
established part of the Margulis conjecture on ergodic actions on the
homogeneous space $SL_k(\R)/SL_k(\Z)$, for $k\ge 3$ (see [9]). 
It was previously well-known
that such a result would have implications to Diophantine questions, 
including to the Littlewood conjecture.  
Their sophisticated proof uses, among others, deep tools from
algebra and from the theory of dynamical systems, involving in
particular the important work developed by Ratner (see for instance [14]).

\medskip

The main purpose of the present Note is to provide a new, short and elementary,
positive answer to a strong form of Question 1. We
will only make use of the basic theory of continued fractions. 
Furthermore, our approach is constructive and allows us to give, 
for any real number $\alpha$ in $\Bad$, continuum many 
explicit examples of pairs 
$(\alpha, \beta)$ of numbers in
$\Bad$ satisfying the Littlewood conjecture, 
with $1$, $\alpha$ and $\beta$ linearly independent over 
the rationals.

\vskip 8mm

\centerline{\bf 2. Main results}

\vskip 6mm

Before stating our main result, we recall the obvious fact that, for any
given $\alpha$ and $\beta$ in $\Bad$, there exists a positive constant
$c(\alpha, \beta)$ such that
$$
q \cdot \Vert q \alpha \Vert \cdot \Vert q \beta \Vert \ge 
{c(\alpha, \beta) \over q }, \eqno (2.1)
$$
for any positive integer $q$. 
Our Theorem 1 gives a positive answer to a strong
form of Question 1 and solves a question posed by de Mathan at the end of
[10].

\proclaim Theorem 1.
Let $\varphi$ be a positive, non-increasing function defined on the set of
positive integers and satisfying $\varphi(1) = 1$, 
$\lim_{q \to + \infty} \, \varphi(q) = 0$ and 
$\lim_{q \to + \infty} \, q\varphi(q) = +\infty$.
Given $\alpha$ in $\Bad$, there exists an uncountable subset
$B_{\varphi}(\alpha)$ of $\Bad$ such that,
for any $\beta$ in $B_{\varphi}(\alpha)$, there exist infinitely many
positive integers $q$ with
$$
q \cdot \Vert q \alpha \Vert \cdot \Vert q \beta \Vert \le 
{1 \over q \cdot \varphi(q)}.
\eqno (2.2)
$$
In particular, the Littlewood conjecture holds for 
the pair $(\alpha, \beta)$ for any $\beta$ in $B_{\varphi}(\alpha)$. 
Furthermore, the set $B_{\varphi}(\alpha)$ can be effectively
constructed.

To the best of our knowledge,
the first explicit examples of independent pairs of real numbers
$(\alpha, \beta)$ satisfying the Littlewood conjecture, with $\alpha$
and $\beta$ both lying in $\Bad$, have been 
recently given by de Mathan in [10]. In particular, for any quadratic
real number $\alpha$, the method introduced by de Mathan allows him to
construct an independent $\beta$ in $\Bad$ such that the pair
$(\alpha,\beta)$ satisfies the Littlewood conjecture. However, his results
yield a positive answer to Question 1 only for a very restricted class
of real numbers $\alpha$.

The proof of Theorem 1 is elementary,
in the sense that it rests only on the theory of
continued fractions. For given $\alpha$ and $\varphi$, we construct
inductively the sequence of partial quotients of a suitable real
number $\beta$ such that (2.2) holds for the  
pair $(\alpha, \beta)$. This sequence can easily be explicited, as we 
show now.

Throughout this Note, we identify any finite or infinite word
$W = w_1 w_2 \ldots$ on the alphabet $\{1, 2, \ldots\} = \Z_{\ge 1}$ 
with the sequence of partial quotients $w_1, w_2, \ldots $
Further, if $U = u_1 \ldots u_m$ and 
$V = v_1 v_2 \ldots$ are words on $\Z_{\ge 1}$, with $V$ finite or infinite,
then $[0; U, V]$ denotes the continued fraction 
$[0; u_1, \ldots, u_m, v_1, v_2,\ldots]$. 
The mirror image of any finite word $W = w_1 \ldots w_m$
is denoted by $\oW := w_m \ldots w_1$.

\proclaim Theorem 2. Let $M \ge 2$ be an integer
and $\eps$ be a positive real number with $\eps < 1$.
Let $\alpha := [0; a_1, a_2, \ldots]$   
be in $\Bad$ with partial quotients
bounded from above by $M$. For any positive integer $n$, denote
by $A_n$ the finite word $a_1a_2\ldots a_n$. 
Let $(t_i)_{i \ge 1}$ be any sequence with values in the set 
$\{M+1,M+2\}$, and let
$(n_i)_{i \ge 1}$ be any sequence of positive integers satisfying
$$
\liminf_{i \to + \infty} \, {n_{i+1} \over n_i} 
> { 4 \log (M+3) \over \eps\log 2}. \eqno (2.3)
$$
Set 
$$
\beta = [0; \oA_{n_1}, t_1, \oA_{n_2}, t_2,
\oA_{n_3}, t_3, \ldots].  \eqno (2.4)
$$
Then, $1$, $\alpha$ and $\beta$ are linearly independent over 
the rationals, and there exist 
infinitely many positive integers $q$ such that
$$
q \cdot \Vert q \alpha \Vert \cdot \Vert q \beta \Vert 
\le {1 \over q^{1 - \eps}}.  \eqno (2.5)
$$
In particular, the Littlewood conjecture holds for the pair 
$(\alpha, \beta)$.

With a slight change in their construction, we may ensure that the
real numbers $\beta$ satisfying the conclusion of Theorem 2 are
transcendental. Indeed, keep the notation of that theorem
and set $B_1 = \oA_{n_1}$ and 
$B_j := \oA_{n_1} t_1 \oB_1 \oA_{n_2} t_2 \oB_2 \ldots
\oB_{j-1} \oA_{n_j}$ for any $j \ge 2$. Then, the real number
$$
\beta = [0; \oA_{n_1}, t_1, \oB_1, \oA_{n_2}, t_2, \oB_2
\oA_{n_3}, t_3, \oB_3 \ldots]
$$
begins in infinitely many palindromes, hence, by Theorem 1
from [1], it is transcendental. To reach the full
conclusion of Theorem 2 with these $\beta$, it is then sufficient to
slightly weaken (2.3).

\medskip

We point out that
a slight modification of the proof of Theorem 1 yields
the following result.

\proclaim Theorem 3.
Let $\varphi$ be as in Theorem 1. Let $M$ be a positive real number.
Let ${\cal A}$ be a countable subset of $\Bad$ such that the partial
quotients of every element of ${\cal A}$ are bounded by $M$.
Then, there exists an uncountable subset  
$B_{\varphi}({\cal A})$ of $\Bad$ such that,  
for any $\alpha$ in ${\cal A}$ and
for any $\beta$ in $B_{\varphi}({\cal A})$, there exist infinitely many
positive integers $q$ with
$$
q \cdot \Vert q \alpha \Vert \cdot \Vert q \beta \Vert \le 
{1 \over q \cdot \varphi(q)}.
$$
Furthermore, the set $B_{\varphi}({\cal A})$ can be effectively constructed.

\noindent
To establish Theorem 3, it is sufficient to follow the proof of
Theorem 1, but, instead of working with the same $\alpha$ at each step,
to work alternatively with each element of ${\cal A}$.  
We omit the details. 

\medskip

Actually, the method for proving Theorem 1 
gives us much freedom, and allows us to
get various results in the same spirit as Theorem 2. Some of them will
be stated in Section 4, with a particular focus on the case when
$\alpha$ and $\beta$ are equivalent real numbers. Section 5 is devoted
to additional remarks and comments.

\vskip 8mm

\centerline{\bf 3. Proofs of Theorems 1 and 2}

\vskip 6mm

For the reader convenience, we recall some well-known results from
the theory of continued fractions, 
whose proofs can be found e.g. in the book of Perron [11].

\medskip

\proclaim Lemma 1.
Let $\alpha = [0; a_1, a_2, \ldots]$ be a real number with convergents
$(p_j/q_j)_{j \ge 1}$. Then, for any $j \ge 2$, we have
$$
{q_{j-1} \over q_j } = [0; a_j, a_{j-1}, \ldots , a_1].
$$

\medskip

\proclaim Lemma 2.
Let $\alpha = [0; a_1, a_2, \ldots]$ and $\beta =
[0; b_1, b_2, \ldots]$ be real numbers. 
Assume that there exists a positive integer $n$ such that
$a_i = b_i$ for any $i=1, \ldots, n$. We then have
$|\alpha - \beta| \le q_n^{-2}$, where $q_n$ denotes the denominator
of the $n$-th convergent to $\alpha$.

\medskip

For positive integers $a_1, \ldots, a_m$, denote 
by $K_m (a_1, \ldots, a_m)$ the denominator of the rational number
$[0; a_1, \ldots, a_m]$. It is commonly called a {\it continuant}.

\medskip

\proclaim Lemma 3.
For any positive integers $a_1, \ldots, a_m$ and any integer $k$ with
$1 \le k \le m-1$, we have
$$
K_m (a_1, \ldots , a_m) = K_m (a_m, \ldots, a_1)
$$
and
$$
\eqalign{
K_k (a_1, \ldots, a_k) \cdot K_{m-k} (a_{k+1}, \ldots, a_m)
& \le K_m (a_1, \ldots , a_m) \cr
& \le 2 \, K_k (a_1, \ldots, a_k) \cdot K_{m-k} (a_{k+1}, \ldots, a_m). \cr}
$$

\medskip

\proclaim Lemma 4.
Let $(a_i)_{i \ge 1}$ be a sequence of positive integers at most equal
to $M$. For any positive integer $n$, we have
$$
2^{(n-1)/2} \le K_n (a_1, \ldots, a_n) \le (M+1)^n.
$$

\medskip

We further need the following auxiliary result.  

\medskip

\proclaim Lemma 5 . Let $M$ be a positive real number. Let 
$\alpha = [0; a_1, a_2, \ldots]$ and $\beta =
[0; b_1, b_2, \ldots]$ be real numbers whose partial quotients are
at most equal to $M$. 
Assume that there exists a positive integer $n$ such that  
$a_i = b_i$ for any $i=1, \ldots, n$ and $a_{n+1}\not=b_{n+1}$. 
Then, we have
$$
|\alpha - \beta| \ge {1\over (M+2)^3 q_n^2},
$$ 
where $q_n$ denotes the denominator
of the $n$-th convergent to $\alpha$.

\medskip

\noindent {\bf Proof.} Set $\alpha'=[a_{n+1}; a_{n+2},\ldots]$ and 
$\beta'=[b_{n+1}; b_{n+2},\ldots]$. Since $a_{n+1}\not= b_{n+1}$, we have 
$$
\vert\alpha'-\beta'\vert\ge 1 - [0; 1, M+1] = 
{1 \over M+2}. \eqno (3.1)  
$$
Furthermore, since the partial quotients of both $\alpha$ and $\beta$ are
bounded by $M$, we immediately obtain   
$$
\alpha'\le M+1 \quad{\rm{and}} \quad \beta'\le M+1. \eqno (3.2)
$$
Denote by $(p_j / q_j)_{j \ge 1}$ 
the sequence of convergents to $\alpha$. Then, the theory of continued
fractions gives that 
$$
\alpha={p_n\alpha'+p_{n-1}\over q_n\alpha'+q_{n-1}} \quad {\rm and} \quad
\beta={p_n\beta'+p_{n-1}\over q_n\beta'+q_{n-1}},
$$
since the first $n$-th partial quotients of $\alpha$ and
$\beta$ are assumed to be the same.  
We thus obtain 
$$
\vert\alpha-\beta\vert=
\left\vert{p_n\alpha'+p_{n-1}\over q_n\alpha'+q_{n-1}}-
{p_n\beta'+p_{n-1}\over q_n\beta'+q_{n-1}}\right\vert
= \left\vert{ \alpha'-\beta'
\over (q_n\alpha'+q_{n-1})(q_n\beta'+q_{n-1})}\right\vert\cdot
$$
Together with (3.1) and (3.2), this yields  
$$\vert\alpha-\beta\vert\ge {1\over (M+2)^3 q_n^2},$$
concluding the proof of the lemma. 
\cqfd

\medskip

We can now proceed with the proofs of Theorems 1 and 2.

\medskip

\noindent {\bf Proof of Theorem 1.}
Write $\alpha =[0; a_1, a_2, \ldots, a_k, \ldots]$.
We first construct inductively a rapidly increasing sequence
$(n_j)_{j \ge 1}$ of positive integers.  
We set $n_1 = 1$ and we proceed with the inductive step.
Assume that $j \ge 2$ is such that $n_1, \ldots , n_{j-1}$ have been
constructed. Then, we choose $n_j$ sufficiently large in order that  
$$\varphi(2^{(m_j-1)/2})\le {1\over 4}\cdot
\left({1\over (M+3)^{m_{j-1}+1}}\right)^2, \eqno (3.3)$$
where $m_j=n_1+n_2+\ldots +n_j+(j-1)$. Such a choice  
is always possible since $\varphi$ tends to zero at infinity and since
the right hand side of (3.3) depends only on $n_1,n_2,\ldots,n_{j-1}$. 

Our sequence $(n_j)_{j \ge 1}$ being now constructed, 
for an arbitrary integer sequence  ${\bf t}=(t_k)_{k\ge 1}$ 
with values in $\{M+1,M+2\}$, we set 
$$
\eqalign{
\beta_{\bf t} = & [0; b_1, b_2, \ldots] \cr
= & [0; a_{n_1}, \ldots, a_1, t_1, 
a_{n_2}, \ldots, a_1, t_2, 
a_{n_3}, \ldots, a_1, \ldots , a_1, t_{j-1}, a_{n_j}, \ldots]. \cr}  
$$
Then, we introduce the set 
$$
B_{\varphi}(\alpha)=
\left\{\beta_{\bf t}, \; {\bf t}\in\{M+1,M+2\}^{\Z_{\ge 1}}\right\}.  
$$
Clearly, the set $B_{\varphi}(\alpha)$ is uncountable.

Let $\beta$ be in $B_{\varphi}(\alpha)$. It remains for us to prove 
that (2.2) with this pair $(\alpha, \beta)$ holds for 
infinitely many integers $q$.
Denote by $(p_j / q_j)_{j \ge 1}$ (resp. by $(r_j / s_j)_{j \ge 1}$)
the sequence of convergents to $\alpha$ (resp. to $\beta$). 
We infer from Lemma 1 that
$$
{s_{m_j - 1} \over s_{m_j} } = [0; a_1, \ldots , a_{n_j},
t_{j-1}, a_1, \ldots, a_{n_{j-1}}, t_{j-2},\ldots, t_1, a_1,
\ldots, a_{n_1}], \eqno(3.4)
$$
which, by Lemma 2, yields
$$
\Vert s_{m_j} \alpha \Vert \le s_{m_j} \, q_{n_j}^{-2}. \eqno (3.5)
$$
Now, we proceed to bound
$s_{m_j} \, q_{n_j}^{-2}$ from above.

By Lemma 4, we have
$$
s_{m_j} \ge 2^{(m_j-1)/2}
$$ 
and, since the partial quotients of $\beta$ are bounded by $M+2$, 
we also get
$$
K_{m_{j-1} + 1} (b_1, \ldots, b_{m_{j-1} + 1})<(M+3)^{m_{j-1} + 1}. 
$$
Thus, using that $\varphi$ is non-increasing,
inequality (3.3) implies that
$$
4 \, \varphi (s_{m_j}) 
\le  K_{m_{j-1} + 1} (b_1, \ldots, b_{m_{j-1} + 1})^{-2}
= K_{m_j - n_j} (b_1, \ldots, b_{m_j - n_j})^{-2} \eqno (3.6)
$$
holds. However, we infer from Lemma 3 that
$$
s_{m_j} \le 2 \, K_{m_j - n_j} (b_1, \ldots, b_{m_j - n_j}) \,
K_{n_j} (b_{m_j - n_j + 1}, \ldots, b_{m_j})  \eqno (3.7)
$$
and
$$
K_{n_j} (b_{m_j - n_j + 1}, \ldots, b_{m_j}) 
= K_{n_j} (a_1, \ldots, a_{n_j}) = q_{n_j}. \eqno (3.8)
$$
Combining (3.6), (3.7) and (3.8), we obtain that
$$
q_{n_j}^{-2} \le {1 \over s_{m_j}^2 \, \varphi (s_{m_j})},
$$
which, together with (3.5), yields
$$
s_{m_j} \cdot \Vert s_{m_j} \alpha \Vert 
\cdot \Vert s_{m_j} \beta \Vert 
\le \Vert s_{m_j} \alpha \Vert \le s_{m_j} \, q_{n_j}^{-2}
\le {1 \over s_{m_j} \, \varphi(s_{m_j})}. 
$$
This shows that (2.2) holds for infinitely many
positive integers $q$ and completes the proof of Theorem 1. \cqfd

\medskip

We now turn to the proof of Theorem 2.

\medskip

\noindent {\bf Proof of Theorem 2.}
With the notation of Theorem 1, we have $\varphi (q) = q^{-\eps}$ for any
positive integer $q$, thus, inequality (3.6) becomes
$$
s_{m_j} \ge \bigl( 2 \, K_{m_{j-1} + 1}
(b_1, \ldots, b_{m_{j-1} + 1})  \bigr)^{2/ \eps}.  \eqno (3.9)
$$
To satisfy (3.9), it follows from 
(3.7), (3.8) and the equality $m_j=m_{j-1}+n_j+1$ that 
it is sufficient to choose $n_j$ such that
$$
2^{(m_{j-1} + n_j)/2} \ge \bigl( 2 (M+3)^{m_{j-1} + 1} \bigr)^{2/\eps},
$$ 
and thus, such that
$$
n_j \ge m_{j-1}\left({4\log (M+3)\over \eps\log{2}}-1\right) +
{4\log(2(M+3))\over \eps\log{2} }. \eqno (3.10)
$$
Our assumption (2.3) implies that (3.10) is satisfied for any
sufficiently large $j$. Consequently, (2.5) holds 
with $\beta$ given by (2.4) 
for any integer $q = s_{m_j}$ large enough. 

It thus only remains to prove that 
$1$, $\alpha$ and $\beta$ are independent 
over the rationals. Therefore, we assume that they are dependent and 
we aim at deriving a contradiction. Let
$(A,B,C)$ be a non-zero integer triple satisfying 
$$
A\alpha+B\beta+C=0. 
$$ 
Then, for any positive integer $q$, we have
$$
\Vert qA\alpha\Vert=\Vert qB\beta\Vert.
$$
In particular, we get
$$
\Vert s_{m_j}A\alpha\Vert= \Vert s_{m_j}B\beta\Vert\leq \vert
B\vert\cdot 
 \Vert s_{m_j}\beta\Vert\ll {1 \over s_{m_j}}, \eqno (3.11)
 $$
 for any $j \ge 2$. Here and below, the constant implied by $\ll$ does
 not depend on $j$.

On the other hand, we have constructed the sequence $(n_j)_{j\ge 1}$ 
in order to guarantee that  
$$
\vert s_{m_j}\alpha-s_{m_j-1}\vert\le  
{1\over s_{m_j}\varphi(s_{m_j})}.\eqno (3.12)
$$
Since by assumption $b_{m_{j-1}+1}=t_{j-1}$ lies in the set $\{M+1,M+2\}$, 
we have $b_{m_{j-1}+1}\not=a_{n_j+1}$. 
Then, (3.4) and Lemma 5 imply that 
$$
\vert s_{m_j}\alpha-s_{m_j-1}\vert\ge
{s_{m_j}\over (M+5)^3 q_{n_j}^2}. 
$$
Moreover, by Lemmas 3 and 4, we obtain 
$$
\eqalign{
s_{m_j}=& \ens K_{m_j}(a_1,\ldots, a_{n_j},b_{m_{j-1} + 1}, a_1, \ldots,
a_{n_{j-1}}, 
b_{m_{j-2} + 1}, a_1,\ldots, a_{n_1})\cr
\ge & \ens K_{n_j}(a_1,\ldots, a_{n_j})\cdot
K_{m_j-n_j}(b_{m_{j-1} + 1}, a_1, \ldots,
a_{n_{j-1}}, 
b_{m_{j-2} + 1}, a_1,\ldots, a_{n_1})\cr
\ge & \ens q_{n_j}2^{m_{j-1}/ 2},}$$ 
hence, we get 
$$
\vert s_{m_j}\alpha-s_{m_j-1}\vert\gg
{2^{m_{j-1}}\over s_{m_j}}. \eqno (3.13)
$$
For $j$ large enough, we deduce from (3.12) that 
$$
\vert s_{m_j}A\alpha- s_{m_j-1} A \vert<{1\over2},
$$  
thus,
$$
\Vert s_{m_j}A\alpha\Vert=\vert s_{m_j}A\alpha- s_{m_j-1} A\vert=
\vert A\vert\cdot\vert s_{m_j}\alpha-s_{m_j-1}\vert.
$$
By (3.13), this yields
$$
\Vert s_{m_j}A\alpha\Vert\gg {2^{m_{j-1}}\over s_{m_j}},
$$
which contradicts (3.11). 
This concludes the proof of Theorem 2. 
\cqfd

\vskip 8mm
\goodbreak

\centerline{\bf 4. Pairs of equivalent numbers}

\vskip 6mm

Two real irrational  
numbers $\alpha$ and $\beta$ are said to be equivalent
(resp. equal up to a rational homography)
if there exist integers $a$, $b$, $c$ and $d$ 
with $\vert ad-bc \vert=1$ (resp. with $\vert ad - bc \vert \not= 0$)
such that 
$$
\beta={a\alpha+b\over c\alpha+d}.
$$

A classical result (see e.g. [11]) asserts
that two real numbers are equivalent 
if, and only if, their continued fraction expansions coincide, up to
finitely many partial quotients.
Consequently, if $\alpha$ is in $\Bad$, then this is also the case for
any real number $\beta$ 
equivalent to $\alpha$. Note that if $\alpha$
is a quadratic real number and if $\beta$ is a real number 
equivalent to $\alpha$, then $\alpha$ and $\beta$ are dependent over
the rational integers. Thus, the Littlewood conjecture holds 
obviously for
any pair of quadratic, equivalent real numbers. 
Moreover, it is
easy to see that if $\alpha$ denotes a non-quadratic irrational number and
if $\beta$ is equivalent to $\alpha$, 
then $1$, $\alpha$ and $\beta$
are independent over the rationals, except if there exists an integer
$m$ such that $\beta=\pm \alpha + m$. 
In particular, $\alpha$ and $1/\alpha$ are equivalent and independent,
for any non-quadratic irrational number $\alpha$.

In the present section, 
we ask wether the Littlewood conjecture is true for
any pair of equivalent real numbers. 
This seemingly innocuous problem is still open, and nothing more is known
on it than on the general conjecture, up to the following remark:
the Littlewood conjecture is true for the
pair $(\alpha, 1/\alpha)$ as soon as $\alpha$ is well approximable
by quadratic numbers [10] (this observation originates in the work
of M. Queff\'elec [13] where she proved the transcendence
of the Thue--Morse continued fraction). 
Actually, this result can be slightly refined: under the same assumption on 
$\alpha$, the Littlewood conjecture is true for the
pair $(\alpha, \beta)$, where $\beta$ is any number 
equivalent to $\alpha$.
We give an explicit related statement in Theorem 4 below and 
describe in Theorem 5 another class of real numbers $\alpha$ such that
the Littlewood conjecture holds for any pair of 
equivalent parameters $(\alpha, \beta)$.

\medskip

In the sequel, we denote by $|W|$ the length of a finite
word $W$. Furthermore, for any positive rational number
$x$, we denote by $W^x$ the word
$W^{[x]}W'$, where $W'$ is the prefix of
$W$ of length $\left\lceil(x- [x])\vert W\vert\right\rceil$ and $\lceil y
\rceil$ denotes the least integer greater than or equal to $y$.

\proclaim Theorem 4. 
Let $\alpha$ be in $\Bad$ and denote 
by $(p_n/q_n)_{n \ge 1}$ the sequence of its convergents.
Assume that there exist a positive rational number $x$
and a sequence of finite words $(U_k)_{k \ge 1}$
such that, for every $k \ge 1$, the continued fraction expansion of
$\alpha$ begins in $[0; U_k, U_k^x]$ and $|U_{k+1}| > |U_k|$.
Set further $M = \limsup_{\ell \to + \infty} \, q_{\ell}^{1/\ell}$ and
$m = \liminf_{\ell \to + \infty} \, q_{\ell}^{1/\ell}$.
If we have
$$
x \ge 1 \quad {\it or} \quad
x > {1 \over 2} \cdot {\log M \over \log m},  \eqno (4.1)
$$
then the Littlewood conjecture is true for the pair $(\alpha, \beta)$,
where $\beta$ is any number equal to $\alpha$ up to a rational homography.

\noindent{\bf Proof.} We content ourselves to
outline the proof. We first recall the dual form of the Littlewood
conjecture (see Lemma 5 from [3]). Given two real numbers 
$\alpha$ and $\beta$, then (1.2) is equivalent to the following
equality 
$$
\inf_{(A,B)\in\Z\times\Z \setminus\{(0,0)\}} \, 
\max\{\vert A\vert,1\}\cdot\max\{\vert B\vert,1\}\cdot \Vert A \alpha 
+ B\beta\Vert = 0.
$$ 

Let $\alpha$ be in $\Bad$ and denote 
by $(p_n/q_n)_{n \ge 1}$ the sequence of its convergents.
For any positive integer $k$, the quadratic number
$\alpha_k := [0; U_k, U_k, \ldots]$ is very close to $\alpha$.
Setting $r_k := |U_k|$, we get
$q_{r_k - 1} \alpha_k^2 + (q_{r_k} - p_{r_k - 1}) \alpha_k
- p_{r_k} = 0$ and
$$
|q_{r_k - 1} \alpha^2 + (q_{r_k} - p_{r_k - 1}) \alpha
- p_{r_k}| \ll q_{r_k} |\alpha_k - \alpha| \ll 
q_{r_k} \, q_{(1+x)r_k}^{-2},
$$
where, as below, the numerical constant implied in $\ll$ depends 
on $\alpha$, but is independent from $k$.  
Then, by (4.1) and Lemma 3, there exists $\eps > 0$ such that
$$
|q_{r_k - 1} \alpha^2 + (q_{r_k} - p_{r_k - 1}) \alpha
- p_{r_k}| \ll q_{r_k} |\alpha_k - \alpha| \ll
q_{r_k}^{-2-\eps}. \eqno (4.2)
$$
Let $\beta = (a \alpha + b)/(c \alpha + d)$ be a number 
equal to $\alpha$ up to a rational homography.
Set $\delta = ad - bc$,
$$
A_k=q_{r_k - 1} \, ,\;\; 
B_k=\delta \bigl( d \bigl((q_{r_k} - p_{r_k - 1})c - q_{r_k - 1} d \bigr)
 +  p_{r_k} c^2 \bigr) 
$$
and
$$ C_k=\delta b  \bigl((q_{r_k} - p_{r_k - 1})c - q_{r_k - 1} d \bigr)
+ \delta p_{r_k} c a.$$
Then, an easy calculation shows that
$$
\Vert A_k \alpha + B_k \beta \Vert
=\left\vert A_k \alpha + B_k \beta - C_k\right\vert   
= {c \over c \alpha + d} \, 
|q_{r_k - 1} \alpha^2 + \ens (q_{r_k} - p_{r_k - 1}) \alpha
- p_{r_k}|. 
$$
Since $\vert A_k\vert\ll q_{r_k}$ and $\vert B_k\vert\ll q_{r_k}$, it 
thus follows from (4.2) that 
$$ 
\max \{\vert A_k\vert,1\} \cdot\max\{\vert B_k\vert,1\} 
\cdot \Vert A_k \alpha + B_k
 \beta \Vert\ll
q_{r_k}^{-\eps}.
$$
This proves that the dual form of the Littlewood conjecture, 
and thus the Littlewood conjecture, holds
for the pair $(\alpha, \beta)$.  \cqfd

\bigskip

In Theorem 4, we used repetition to construct suitable real numbers $\alpha$.
Another useful combinatorial tool is palindromy. We recall that a
palindrome is a finite word $W$ such that $\oW=W$.

\proclaim Theorem 5. 
Let $\alpha$ be in $\Bad$ and denote 
by $(p_n/q_n)_{n \ge 1}$ the sequence of its convergents.
Assume that there exist a positive rational number $x$
and two sequences of finite words $(U_k)_{k \ge 1}$ and $(V_k)_{k \ge 1}$
such that, for every $k \ge 1$, 
the continued fraction expansion of
$\alpha$ begins in $[0; V_k, U_k, \oU_k]$ and $|U_{k+1}| > |U_k| \ge x |V_k|$.
Set further $M = \limsup_{\ell \to + \infty} \, q_{\ell}^{1/\ell}$ and
$m = \liminf_{\ell \to + \infty} \, q_{\ell}^{1/\ell}$.
If we have
$$
x > {3 \over 2} \cdot {\log M \over \log m} - 
{1 \over 2}, \eqno (4.3)
$$
then the Littlewood conjecture is true for the pair $(\alpha, \beta)$,
where $\beta$ is any number equal to $\alpha$ up to a rational homography.

\noindent{\bf Proof.} Let $\beta = (a \alpha + b)/(c \alpha + d)$ be a number 
equal to $\alpha$ up to a rational homography.
Let $k \ge 1$ be an integer and let $P_k/Q_k$ be  
the last convergent to the rational number
$$
{P'_k \over Q'_k} := [0; V_k, U_k, \oU_k, \oV_k].
$$
It follows from Lemma 1 that $P'_k = Q_k$.
Setting $r_k = |U_k|$ and $s_k=\vert V_k\vert$, 
we infer from Lemma 2 that
$$
\max \{ \Vert Q_k \alpha \Vert,
\Vert Q'_k \alpha \Vert \}
\ll \, Q_k \,
q_{s_k + 2 r_k}^{-2}. \eqno (4.4)
$$
Here and below, the constants implied by $\ll$ may depend on $\alpha$ 
and $\beta$, but not on $k$. 
Observe that
$$
\biggl| \beta - {a P'_k + b Q'_k \over c P'_k + d Q'_k} \biggr|
\ll \biggl| \alpha - {P'_k \over Q'_k} \biggr|.
$$
Thus, setting
$$
R_k := |c P'_k + d Q'_k| = |c Q_k + d Q'_k|,
$$
we get $R_k \ll Q_k$ and 
$$
\max \{ \Vert R_k \alpha \Vert,
\Vert R_k \beta \Vert \} 
\ll \, Q_k \, q_{s_k + 2 r_k}^{-2}, \eqno (4.5)
$$
using (4.4).
Furthermore, Lemma 3 implies that 
$$ 
Q_k \ll K_{2(r_k+s_k)}(V_k U_k\, \overline{U_k} \, \overline{V_k})\ll
q_{s_k} \, q_{s_k + 2 r_k}.
$$
Then, it follows from (4.5) that 
$$
R_k \cdot \Vert R_k \alpha \Vert \cdot \Vert R_k \beta \Vert\ll
q_{s_k}^3 \, q_{s_k + 2 r_k}^{-1}. \eqno (4.6)
$$
In virtue of (4.3), this concludes the proof. 
 \cqfd

\bigskip

Actually, a sharper conclusion
than (1.2) holds for the pairs $(\alpha, \beta)$ satisfying the
conclusion of Theorem 5: there exists a positive real number $\eps < 1$, 
depending on $x$, such that (2.5) holds for infinitely many
positive integers $q$. This can be further refined under the strongest
assumption that $\alpha$ begins in arbitrarily large palindromes.

\proclaim Theorem 6. 
Let $\alpha$ be in $\Bad$ such that its continued fraction expansion
begins in infinitely many palindromes. Let $\beta$ be any real
number equal to $\alpha$ up to a rational homography. 
Then, the Littlewood conjecture is true
for the pair $(\alpha, \beta)$ and, moreover, we have
$$
\liminf_{q \to + \infty} \,
q^2 \cdot \Vert q \alpha \Vert \cdot \Vert q \beta \Vert < + \infty.
$$

The proof of Theorem 6 follows the same lines as that of
Theorem 5: it essentially amounts to setting  
$s_k = 0$ in (4.6). 

For $\alpha$ being as in Theorem 6,
the fact that the Littlewood conjecture is true
for the pair $(\alpha, 1/ \alpha)$ has previously been
noticed by M. Queff\'elec in her talk held at the I.H.P. in June 2004.


\vskip 8mm

\centerline{\bf 5. Concluding remarks}

\vskip 6mm

For the reader convenience, 
we reformulate inequality (2.1). For any
given $\alpha$ and $\beta$, both lying in $\Bad$, 
there exists a positive constant
$c(\alpha, \beta)$ such that
$$
q^2 \cdot \Vert q \alpha \Vert \cdot \Vert q \beta \Vert >
c(\alpha, \beta), 
$$
for any positive integer $q$. 
In view of this
and of Theorem 1, we propose the following problem:

\medskip

{\it Question 2.} Given $\alpha$ in $\Bad$, is there any independent
$\beta$ in $\Bad$ so that 
$$
\liminf_{q \to + \infty} \,
q^2 \cdot \Vert q \alpha \Vert \cdot \Vert q \beta \Vert < + \infty \;
?
\eqno (5.1)
$$

\medskip

\noindent 
We observe that (5.1) holds
when $\alpha$ and $\beta$ are linearly dependent over the rationals,
as follows from (1.1).
Furthermore, Theorem 6 gives a positive answer to Question 2
for a restricted class of real numbers $\alpha$.
Apart from this partial result, we do not know the answer
to Question 2.

\medskip

Let $\K$ be any field, and let $X$ be an indeterminate.
We define a norm $| \cdot |$
on the field $\K ((X^{-1}))$ by setting $|0| = 0$ and,
for any non-zero formal power series 
$F(X) = \sum_{h=-m}^{+ \infty} \, f_h X^{-h}$ with $f_m \not= 0$, 
by setting $|F| = 2^m$. We further write $|| F ||$ to denote the norm
of the fractional part of $F(X)$, that is, of the part of the series
which comprises only the negative powers of $X$.
In analogy with the Littlewood conjecture, we may ask whether,
given $F(X)$ and $G(X)$ in $\K ((X^{-1}))$, we have
$$
\inf_{q \in \K[X] \setminus \{0\}}  
\, |q| \cdot \Vert q F \Vert \cdot \Vert q G \Vert = 0.
$$
A negative answer to this question has been obtained by 
Davenport \& Lewis [5] (see also Baker [2]
for an explicit counter-example) when $\K$ is an infinite field. 
The question is still unsolved when $\K$ is a finite field.
We conclude by pointing out that 
our construction can also be applied to solve the analogue of Question 1 
for formal power series defined over an arbitrary field. 
This will be part of a subsequent work. 

\vskip 5mm
\noindent {\bf Acknowledgements.} We are very grateful to 
Bernard de Mathan for many useful remarks. In particular, the present
version of Theorem 5 and its proof is due to him.

\vskip 12mm \goodbreak

\centerline{\bf References}

\vskip 7mm

\item{[1]}
B. Adamczewski and Y. Bugeaud,
{\it Palindromic continued fractions}. Preprint.

\sm

\item{[2]}
A. Baker,
{\it On an analogue of Littlewood's diophantine approximation
problem}, Michigan Math. J. 11 (1964), 247--250.

\sm

\item{[3]} 
J. W. S. Cassels and H. P. F. Swinnerton-Dyer,
{\it On the product of three homogeneous linear forms and indefinite
ternary quadratic forms}, Philos. Trans. Roy. Soc. London,
Ser. A, 248 (1955), 73--96.

\sm

\item{[4]} 
H. Davenport, P. Erd\H os and W. J. LeVeque,
{\it On Weyl's criterion for uniform distribution}, Michigan
Math. J. 10 
(1963), 311--314.

\sm
\item{[5]}
H. Davenport and D. J. Lewis, 
{\it An analogue of a problem of Littlewood}, Michigan Math. J. 
10 (1963) 157--160.

\sm

\item{[6]}
M. Einsiedler, A. Katok and E. Lindenstrauss,
{\it Invariant measures and the set of exceptions to the
Littlewood conjecture},
Ann. of Math. To appear.

\sm

\item{[7]} 
R. Kaufman, 
{\it Continued fractions and Fourier transforms},
Mathematika 27 (1980), 262--267.

\sm

\item{[8]}
J. E. Littlewood,
Some problems in real and complex analysis.
D. C. Heath and Co. Raytheon Education Co., 
Lexington, Mass., 1968.

\sm

\item{[9]}
G. A. Margulis,
{\it Problems and conjectures in rigidity theory},
Mathematics: frontiers and perspectives, pages
161--174. Amer. Math. Soc., Providence, RI, 2000.

\sm

\item{[10]}
B. de Mathan,
{\it Conjecture de Littlewood et r\'ecurrences lin\'eaires}, 
J. Th\'eor. Nombres Bordeaux 13 (2003), 249--266.

\sm

\item{[11]}
O. Perron,
Die Lehre von den Kettenbr\"uchen,
Teubner, Leipzig und Berlin, 1929.

\sm

\item{[12]} 
A. D. Pollington and S. Velani,
{\it On a problem in simultaneous Diophantine approximation:
Littlewood's conjecture},
Acta Math. 185 (2000), 287--306.

\sm

\item{[13]}
M. Queff\'elec,
{\it Transcendance des fractions continues de Thue--Morse},
J. Number Theory 73 (1998), 201--211.

\sm

\item{[14]}
M. Ratner,
{\it Interactions between ergodic theory, Lie groups, and number
theory}, Proceedings of the International Congress of Mathematicians,
Vol. 1, 2 (Z\"urich 1994), 157--182, Birkh\"auser, Basel, 1995.

\vskip21mm

\noindent Boris Adamczewski   \hfill{Yann Bugeaud}

\noindent   CNRS, Institut Camille Jordan  
\hfill{Universit\'e Louis Pasteur}

\noindent   Universit\'e Claude Bernard Lyon 1 
\hfill{U. F. R. de math\'ematiques}

\noindent   B\^at. Braconnier, 21 avenue Claude Bernard
 \hfill{7, rue Ren\'e Descartes}

\noindent   69622 VILLEURBANNE Cedex (FRANCE)   
\hfill{67084 STRASBOURG Cedex (FRANCE)}

\vskip2mm
 
\noindent {\tt Boris.Adamczewski@math.univ-lyon1.fr}
\hfill{{\tt bugeaud@math.u-strasbg.fr}}

\bye